\input amstex
\documentstyle{amsppt}
\magnification=\magstep0
\define\cc{\Bbb C}

\define\N{\Bbb N}
\define\jj{\Bbb J}

\define\A{\Cal A}
\define\h{\Cal D}

\define\m{\Cal M}

\define\f{\Cal S}

\define\n{\Cal F}

\define\la{\lambda}
\define\om{\omega}

\define\e{\varepsilon}
\define\va{\varphi }
\define\CB#1{\Cal C_b(#1)}
\define\st{\subset }
\define\al{\alpha}
 \topmatter
 \title
 Existence   of bounded uniformly continuous mild solutions on $\Bbb{R}$ of evolution equations and their asymptotic behaviour.
\endtitle
 \author
 Bolis Basit and Hans G\"{u}nzler
\endauthor
\abstract{ We prove that   $u'= A u +  \phi  $ has on $\Bbb{R}$ a mild solution $u_{\phi}\in BUC (\Bbb{R},X)$ (that is bounded  and uniformly continuous), where $A$ is the generator of a  $C_0$-semigroup
   on  the Banach space ${X}$ with resolvent satisfying  $||R(it,A)||= O(|t|^{-\theta})$,  $|t|\to \infty $, with some $\theta > \frac{1}{2}$, $\phi\in L^{\infty} (\Bbb{R},{X})$ and  $i\,sp (\phi)\cap \sigma (A)=\emptyset$.   As a consequence it is shown that if ${\Cal F}$ is the  space of almost periodic, almost automorphic,
          bounded Levitan almost periodic or certain classes of  recurrent functions  and $\phi$ as above is such that $M_h \phi:=(1/h)\int_0^h \phi (\cdot+s)\, ds \in \Cal {F}$  for each $h >0$, then $u_{\phi}\in \Cal {F}\cap BUC (\Bbb{R},X)$. These results seem new  and strengthen several recent theorems.}
 \endabstract
\endtopmatter
\rightheadtext{Existence of bounded  mild solutions} \leftheadtext{ Basit and
 G\"unzler}
 \TagsOnRight

\document
\baselineskip=22pt

\head{\S 1.  Introduction}\endhead

 In the following \footnote {AMS subject classification 2010: \,\,\,\, Primary  {47D06, 43A60}\,\,\,\, Secondary {43A99,34B27, 35B34, 35C15, 47A10}.
 \newline\indent Key words:  Evolution equations, non-resonance, global solutions, bounded
             uniformly continuous solutions, almost periodic, almost automorphic, Green function.} a linear translation invariant subspace $\n$ of $L^{\infty} ({\Bbb R},X)$
 with complex Banach space $X$ and linear $A: D(A)\to X$ will be called admissible for

(1.1) \qquad $u'= A u +  \phi  $ on   $  {\Bbb R}$,

 \noindent if for every  $\phi\in \n$ with  ($sp=$  Beurling spectrum)

(1.2) \qquad  $i\,sp(\phi)\cap \sigma (A) =\emptyset$,

\noindent (1.1) has a mild solution $u_{\phi}\in \n $ (see (3.2)).
The definitions vary, see    [16, p. 126], [15, p. 167], [21, Definition 11.3, p. 287, 306], [20, p. 401],  [17, p. 248].   A very good survey of previous results here can be found in the introduction  in Phong-Sch\"{u}ler [20].

In [4, Theorem 6.5 (ii)]  with results on the operator equation  $AX-XB=C$ of
[19] it has been shown that $BUC ({\Bbb R},X) = \{f :{\Bbb R}\to X$  bounded uniformly
continuous$\}$ is admissible  if  $A$ is the generator of a holomorphic $C_0$-semigroup $T$ with
sup $_{t > 0} ||T(t)|| < \infty$.

In [20] it is shown that one has equivalence   between admissibility of a translation invariant subspace $\n \st BUC ({\Bbb R},X)$ with respect to (1.1) and the unique solvability of a special operator equation of Lyapunov's type $AX-X\h_{\n}=-\delta_0$, where $\h_{\n}$ is the restriction of the operator $\h:=\frac{d}{dt}$ to $\n$ and $\delta_0 \phi =\phi(0)$.

 Using  this  and spectral properties  of the sum of commuting operators from [1,
Theorem 7.3] a new approach to admissibility has been given in [19], [20] and [17, results
in section 3] for $\n\st  BUC({\Bbb R},X) $  if either $sp (f)$ is compact for $f\in\n$ or $T$ holomorphic or
       $T$ admits exponential dichotomy.

In Theorem 5.2 below  we establish the existence of a bounded uniformly continuous
mild solution  $u_{\phi}$ on ${\Bbb R}$ of the form  $u_{\phi} = G * \phi$  with  $G  \in  L^1({\Bbb R},L(X))$ for any
$\phi \in L^{\infty} ({\Bbb R},X)$ with (1.2), when $A$ is the generator of a $C_0$-semigroup
$T$ with resolvent satisfying $||R(it,A)|| = O(|t|^{-\theta})$, $|t|\to \infty$, with some $\theta > \frac{1}{2}$; $T$ holomorphic is  the case $\theta =1$ . So if additionally  sup$_{t>0} ||T(t)|| < \infty$, for each $x\in  X$
the unique mild solution of the Cauchy problem $u(0) = x$  on $[0,\infty)$ is $\in
BUC ({\Bbb R}_+,X)$. Comparing Theorem 5.2 with the Non-resonance Theorem 5.6.5 of
[2], our result is a 3-fold extension: It gives solutions on ${\Bbb R}$ instead of $[0,\infty)$,
$T$ need not be bounded, and T need not be holomorphic (special case
$\theta = 1$). Theorem 5.6.5 of [2] is more general since it uses the (smaller)
half-line spectrum instead of the Beurling spectrum used in (1.2); however
in the important cases of almost periodic, almost automorphic, Levitan
almost periodic or recurrent functions these two spectra coincide by  [8, Example 3.8], [9, Corollary  5.2].
Also, the proof of this Non-resonance  Theorem of [2] seems not to be
        extendable to general function classes as in our section 6:

In Theorem 6.3 it is  shown that for any linear BUC-invariant  $\n \st L^1_{loc}
 ({\Bbb R},X)$ closed under uniform convergence and satisfying (6.5) the $L^{\infty}  \cap  \m\n$  is   admissible
for (1.1), and additionally the solution  $u_{\phi} \in \n \cap BUC ({\Bbb R},X)$, with $A$ as in
Theorem 5.2, and  $\m\n =$  first mean extension of $\n$ of (6.1) below.

Examples are $ \n =$ almost periodic functions $ AP=AP ({\Bbb R},X)$,  $\st $ Stepanoff almost
periodic functions  $S^p AP \st  \m AP$, $1 \le p < \infty$ [5, (3.8)], so for bounded
$S^1$-almost periodic $\phi$ with (1.2) the solution $u_{\phi}$ is even Bohr almost
periodic. This generalizes for example  results in [19], [4], [20].
Or  $\n = $ Veech almost automorphic functions $VAA=VAA ({\Bbb R},X)$ [23], $\st  \m VAA$;
then for  $\phi \in  L^{\infty} \cap \m VAA$  the $u_{\phi}  \in VAA \cap BUC$,  and so even
Bochner almost automorphic $ \in  AA $  [26, p. 66], [6, (3.3)]. This generalizes a result of [11, Theorem
4.5] in several directions: The semigroup $T$ need not be holomorphic, $\phi \in \m VAA \cap L^{\infty}$ suffices instead of  $\phi\in  AA$, strictly $\st  VAA \st  \m VAA$
[6,(3.3)], the solution $u_{\phi}$ is in addition  $\in BUC$,  and
in (1.2) the Beurling spectrum can be used instead of the  uniform
spectrum $sp_u$  of [11, p. 3293].
Further examples would be  $\n  =$ bounded Levitan almost periodic functions [6, p. 430,
Proposition 3.4], or linear invariant subspaces  of recurrent
functions [6, p. 427], and various spaces of asymptotic almost periodic
functions (see Examples 6.2).
So also $\n$ with not necessarily compact or uniformly continuous elements
are included; these seem  not to be treatable by the methods used in [19], [4], [20],
 [17], [11].

\head{\S 2.  Notation and Definitions}\endhead

In the following ${\Bbb R}_+ = [0,\infty)$, $X$ is a complex Banach space,  $L(X) =$ the Banach algebra $\{B : X\to X,   B $ linear bounded $\}$ with operator norm
$||B||$, $\h({\Bbb R})$ and  $\f ({\Bbb R})$ contain  Schwartz's complex valued  $C^{\infty }$-functions with compact support respectively  rapidly decreasing derivatives,
$BUC({\Bbb R},X) = \{ f : {\Bbb R}\to X : f$  bounded uniformly  continuous $\}$, $AP = AP({\Bbb R},X)$ almost periodic  functions [2, p. 285],  $VAA = VAA({\Bbb R},X)$ (Veech-) almost automorphic functions [23], $AA = AA({\Bbb R},X)$ Bochner almost automorphic functions  [26, p. 66], [6, p.430], [11, p. 3292]; for
$f \in  L^1_{loc}(J,X)$  $(Pf)(t) : =$
Bochner integral  $\int^t_0 f(s)\, ds$, $\widehat {f}(\lambda) = L^1$-Fourier transform
 $\int_{{\Bbb R}} f(t) e^{- i \lambda\, t}\, dt$,  $f_a(t) =$  translate $f(a+t)$ where defined, $a$ real, $sp
=$  Beurling spectrum (3.4), Proposition 3.3.

\head{\S 3. Preliminaries}\endhead

 We study solutions of the inhomogeneous  abstract evolution equation

(3.1)\qquad $\frac{d u(t)}{dt}= A u(t) +  \phi (t) $,  $t\in  {J}$,

\noindent  where $A: D(A)\to X$ is the generator of a
 $C_0$-semigroup $(T(t))_{t\ge 0}$
on the complex Banach space ${X}$, $J\in\{{\Bbb R}_+,{\Bbb R}\}$ and $\phi\in L^1_{loc} ({J},{X})$.

\noindent By [18, Corollary 2.5, p. 5], it follows that $D(A)$ is dense in $X$ and $A$ is a closed linear operator.

\noindent By a $classical\,\, solution$ of (3.1) we mean a function $u: J\to D(A)$ such that $u \in C^1 (J,X)$ and (3.1) is satisfied.

\noindent  By a $mild\,\, solution$ of (3.1) we mean a $\omega  \in
C(J,X)$  with  $\int _ 0 ^ t  \omega(s)ds  \in  D(A)$ for  $t \in J$ and

(3.2)\qquad $\om(t)  =  \om(0)  + A \int_0^t  \omega(s)\, ds  + \int_0^t \phi(s)ds$, $t  \in J $.

\noindent For $ J = {\Bbb R}_+$  this is the usual definition [2, p. 120].

A classical solution is always a mild solution (see [2, (3.1), p. 110], for
$J = {\Bbb R}_+$, $\phi = 0$).
Conversely, a mild solution with $\omega  \in  C^1(J,X)$ and $\phi \in C(J,X)$ is
a classical solution (as in [2, Proposition 3.1.15]).

In the following we collect some  needed lemmas and propositions.

With translation and the case $J = {\Bbb R}_+ $ [2, Proposition 3.1.16] one can show

\proclaim {Lemma 3.1} If $J \in \{{\Bbb R}_+,{\Bbb R}\}$, $T$, $A$ and $\phi$ are as after (3.1) and

    $(M_h\omega) (t) : = (1/h) \int^h_0 \omega (t+s)\,ds$,

  \noindent   the following 3
    statements are equivalent:

  (i)  $\omega $ is  a mild solution of (3.1),

(ii) $\omega \in C(J,X)$ and for all  $ 0 < h \in{\Bbb R}$   the $M_h \omega$ is a
              classical solution of (3.1),

(iii) for all    $t_0  \in J$  one has

(3.3) \qquad $\om(t)= T(t-t_0)\om(t_0)+\int_{t_0}^t T(t-s)\phi(s)\,ds$,\,\,\,$t\ge t_0$.
\endproclaim

\proclaim{Proposition 3.2}  Let  $F\in L^1({\Bbb R},L(X))$ and  $\phi\in L^{\infty}({\Bbb R},X)$ respectively  $F\in L^1({\Bbb R},X)$ and  $\phi\in L^{\infty}({\Bbb R},\cc)$ respectively  $F\in L^1({\Bbb R},\cc)$ and  $\phi\in L^{\infty}({\Bbb R},X)$. Then

(i) $F*\phi (t):= \int_{{\Bbb R}}F(s)\phi (t-s)\, ds$

\noindent exists as a Bochner integral for $t\in {\Bbb R}$ and
$F*\phi\in BUC({\Bbb R},X)$.

(ii)   If additionally $f \in  L^1({\Bbb R},\cc)$, then the convolution $F*\phi*f $ exists    and is associative.
\endproclaim

\demo{Proof} (i) As in [2, Proposition 1.3.4, p. 24]; uniform continuity of $F*\phi$
     follows by approximating $F$ in $L^1$ by step functions.

(ii) Follows  with the Fubini-Tonelli  theorem (see [14, Satz 3, p. 211]). \P
\enddemo

In the following $sp\,\,$ denotes the $Beurling\, spectrum$, $sp (\phi) = sp_{\{0|R\}} (\phi)$ as defined for example in
 [8, (3.2), (3.3)] case $S=\phi\in L^{\infty} ({\Bbb R},X)$, $V= L^1 ({\Bbb R},\cc)$, $\A= \{0|{\Bbb R}\}$:

(3.4) \qquad $sp_{\{0|{\Bbb R}\}} (\phi):= \{\om\in {\Bbb R}: f\in L^1 ({\Bbb R},\cc), \phi*f =0 $ imply $\widehat {f}(\om)=0\}$.

$sp_B$ is defined in [2, p. 321], $sp_C$ is the Carleman spectrum
[2, p.293/317].

\proclaim {Proposition 3.3}  If  $\phi\in L^{\infty} ({\Bbb R},X)$, then

(3.5) \qquad  $sp (\phi)  =  sp_{\{0|R\}} (\phi)  =  sp_B (\phi)  =  sp_C (\phi) $.

\noindent Moreover, if  $sp (\phi)=\emptyset$, then $\phi = 0$ a.e.
\endproclaim
\noindent See also [8, (3.3), (3.14)].
\demo{Proof} $sp (\phi)  \st  sp_B \phi $:   Assume $\la \in sp (\phi)$.
            To any $\e > 0 $ there is  $h \in L^1({\Bbb R},\cc)$  with $\widehat {h} (\la)\not = 0$ and supp $\widehat {h} \st [\la-\e,\la+\e]$,  we conclude $\phi*h \not =0$. This  implies $\la \in sp_{B} (\phi)$.

 $sp_{B} (\phi)\st sp (\phi)$: Assume  $\la \not\in sp(\phi)$. Then there is $h\in L^1 ({\Bbb R},\cc)$ with $\widehat {h} (\la)\not = 0$ and $\phi*h =0$.  With Wiener's
inversion theorem [10, Proposition 1.1.5 (b), p. 22], there is $h_{\la}\in L^1$ such that $k=h*h_{\la}$ satisfies
 $\widehat{k}=1$ on some neighbourhood $V = (\la-\e,\la+\e)$ of $\la$. Now let $g\in L^1 ({\Bbb R},\cc)$ be such that supp$\, \widehat {g}\st V$. Then $k*g=g$ and $0= \phi*h = \phi*k = \phi *(k*g)= \phi*g$. This implies $\la\not \in sp_{B}(\phi)$ by [2, p. 321].

$ sp_B (\phi) = sp_C(\phi)$:  See  [2, Proposition  4.8.4, p. 321].

 \noindent So, $sp_C (\phi)=\emptyset$ if
$sp (\phi)=\emptyset$.  It follows  that $\phi =$ a.e. by  [2, Proposition  4.8.2 a, p. 319]. \P
\enddemo

\proclaim{Lemma 3.4}  Let  $\phi \in L^{\infty}({\Bbb R},X)$.
If  $F  \in L^1({\Bbb R},\cc)$ or $F\in L^1 ({\Bbb R},L(X))$, then $ sp(F*\phi)  \st   sp (\phi) \cap $ supp $\widehat{F}$.
\endproclaim

\demo{Proof}  $F \in L^1({\Bbb R},L(X))$: to $\lambda \not \in sp  (\phi)$ exists $f \in L^1({\Bbb R},\cc)$ with $\phi*f=0$,
     $\widehat{f}(\lambda)=1$, then  $(F*\phi)*f = F*(\phi*f) = 0$ with Proposition 3.2(ii), so $\lambda \not \in
     sp(F*\phi)$, yielding  $ sp(F*\phi)  \st  sp (\phi) $.

     If  $\lambda \not \in supp \widehat{F}$, there exists $f \in \f({\Bbb R})$ with $(\text {supp\,} \widehat{f}) \cap\text {supp\,} \widehat{F}
     =  \emptyset$ and $\widehat{f}(\lambda)=1$, then $\widehat {F*f} = \widehat {F}\widehat {f} = 0$, and so $F*f = 0$;  this gives
      $ 0 =  (F*f)*\phi = F*(f*\phi) = F*(\phi*f) = (F*\phi)*f$,  and so $\lambda
      \not \in  sp(F*\phi)$ by Proposition 3.2 (ii), yielding   $sp(F*\phi)  \st \text{ supp\,}\widehat {F}$.
        $  F*f  \in  L^1({\Bbb R},L(X))$ and  $\widehat {F*f} = \widehat {F}\widehat {f}$ follow
       with Fubini-Tonelli.

The proof for the case $F  \in L^1({\Bbb R},\cc)$ is similar. \P
\enddemo

\proclaim{Lemma 3.5}  To $K$ compact  $\st  U $ open  $\st {\Bbb R}$   there exist an open $V$ and
    $0 \le  \va \in \h({\Bbb R})$  with  $K  \st  V  \st  \overline {V}$ compact  $\st  U$, $\va = 1$ on
    $V$  and  supp $\va  \st U $.
\endproclaim
\demo{Proof} With a partition of unity, for example [22,  p. 147, Theorem 6.20], $\va =
     \sum^m_{j=1} \psi_j$  there.    \P
\enddemo
\proclaim{Lemma 3.6} To  $\phi\in L^{\infty}({\Bbb R},X)$, $M$ compact  $\st {\Bbb R}$ with  $M  \cap  sp (\phi)
   =\emptyset$  there exists a sequence  $(\Pi_n)_{n\in\N}$ of $X$-valued trigonometric
    polynomials with

(3.6)     sup$_{n\in\N}  ||\Pi_n||_{\infty}   <  \infty$,

(3.7)     $\Pi_n  \to  \phi$   almost everywhere  on ${\Bbb R}$,

(3.8)   All Fourier exponents of the $\Pi_n  $ are in $  {\Bbb R}\setminus M $,  $n \in \N$.
\endproclaim

\demo{Proof} To $\phi$ there exist $\phi_n \in C({\Bbb R},X)$  with  $\int^n_{-n}||\phi (t)-\phi_n (t)||dt <
2^{-n}$, $\phi_n = 0 $ outside $[-n,n]$  and $||\phi_n||_{\infty} \le ||\phi||_{\infty}$, $n\in\N$,
for example via step functions $h_n$ and the cut-off operation $h_n  \cap  ||\phi||_{\infty} $ of
[14, p. 327, (12), (13)].

Extending the $\phi_n$ with period $2n$ to ${\Bbb R}$, with Fejer summation [2,Theorem 4.2.19] there
exist $X$-valued trigonometric polynomials $\psi_n$ with
\,\, $\int ^n_{-n} ||\phi_n (t)- \psi_n(t)|| dt\,\,  < 2^{-n}$ and \,\,  $||\psi_n||_{\infty}\, \le ||\phi||_{\infty} +1$,\,\,\, $n\in\N$. This implies \,\,\, $\int^n_{-n}\, ||\psi_{n+1}(t) - \psi_n (t)||\, dt\, \le$

\noindent  $4\times  2^{-n}$,
$n\in\N $. So for fixed n the Monotone Convergence theorem for $L^1({\Bbb R},{\Bbb R})$
implies  $\sum^{\infty}_{m=1} ||\psi_{m+1}(t) - \psi_m(t)||  <  \infty$  for almost all $t\in [-n,n]$. Therefore  $lim_{m \to {\infty}} \psi_m(t)$ exists $= : \psi(t) \in X$  for almost
all $t \in [-n,n]$. Since  $\psi_m\to\phi$ in $L^1 ([-n,n],X)$, there exists a subsequence
$\psi_{m_k}\to\phi$  almost everywhere in $[-n,n]$, so $\psi = \phi$ almost everywhere, that is  $\psi_m
\to\phi$  almost everywhere in $[-n,n]$;  $n$ being arbitrary, (3.6) and (3.7) follow for $\Pi_n =
\psi_n$.

(3.8) : With Lemma 3.5 and  $\h ({\Bbb R})  \st  \f({\Bbb R})  =  \widehat {\f}({\Bbb R})$  choose  $\psi\in\f ({\Bbb R})$ with

supp $\widehat {\psi}
\st {\Bbb R}\setminus sp(\phi)$  and  $\widehat {\psi} = 1$ on some open $V$ with $M \st V$.

\noindent Since  $\psi*\gamma_{\omega}  = \widehat {\psi}(\omega) \gamma_{\omega}$, $\gamma_{\omega} (t) : =
e^{i \omega t}$,  $\omega\in {\Bbb R}$,  the $\Pi_n : = \psi_n - \psi_n * \psi$  are trigonometric polynomials satisfying  (3.6) and  $\Pi_n\to (\phi -\phi*\psi)$  almost everywhere on ${\Bbb R}$.
Since   $ sp (\phi) \cap $  supp $\widehat {\psi}   =\emptyset$, one has  $sp(\phi*\psi) \st  sp(\phi) \cap $ supp $\widehat {\psi}  =\emptyset$
by Lemma 3.4; so $\phi*\psi  =  0 $ on ${\Bbb R}$ with
 Proposition 3.3, the $\Pi_n$ satisfy (3.7).
With the above $\psi*\gamma_{\omega} =\widehat {\psi}(\omega) \gamma_{\omega}$  and
$\widehat {\psi} = 1$ on $M$ one gets (3.8) for the $\Pi_n $.   \P
\enddemo

\proclaim{Proposition 3.7} If $f,\widehat{f} \in L^1 ({\Bbb R},X)\cap C ({\Bbb R},X)$, then ([2, Theorem 1.8.1 d), p. 45])

$f(t) = \frac{1}{2\pi}\int_{{\Bbb R}} \widehat{f} (\la) e^{i\la t}\, d\la$, $t\in {\Bbb R}$.
\endproclaim

\head{\S 4.  Construction of a Green function}\endhead

In all of the following we assume

(4.1)\qquad  $T$ is a $C_0$-semigroup with generator $A:D(A)\to X $;

the $D(A)$ is then linear and dense in $X$, $A$ is a closed linear operator
[18, p. 1, Corollary 2.6  p. 5].
$\sigma(A) : =$ spectrum of $A$,  $R(\cdot,A) : \cc\setminus \sigma(A) \to L(X)$  resolvent of
$A$  [2, Appendix B, p. 462].

(4.2)\qquad   $K : = (1/i)(i{\Bbb R} \cap \sigma(A))$, $R(t) : = R(it,A)$, $t \in {\Bbb R}\setminus K$.

We assume further with $K$, $R$ of (4.2)

(4.3)\qquad there exists  $a \in (0,\infty)$  with $K  \st  (-a,a) = : I_a $,

(4.4) \qquad  there exists $\theta \in (\frac{1}{2},\infty)$ with   sup $\{|t|^{\theta} ||R(t)|| : t \in {\Bbb R}\setminus I_a\} <                                                        \infty $.

\noindent ``$T$ holomorphic" [2, Definition 3.7.1, p. 152]  is by
        [2, Corollary 3.7.18, p. 160] the special case $\theta = 1$ in (4.4).

For later estimates we need $\theta < 1$ :

\proclaim{Lemma 4.1}  (4.1), (4.2), (4.3), (4.4)  imply

(4.5)  \qquad    there exists $\delta \in (\frac{1}{2},1)$ with $\eta : =$  sup $\{|t|^{\delta} ||R(t)|| : t \in {\Bbb R}\setminus I_a\} <    \infty $.
\endproclaim

\demo{Proof} $R$ is continuous on ${\Bbb R}\setminus I_a$ by [2, Corollary B.3, p. 463], $|t|^{\delta} \le
  |t|^{\theta}$      if $|t| \ge 1$,  $0 < \delta \le \theta$,  $a> 0 $.  \P
\enddemo

\proclaim{Lemma 4.2} (4.1), (4.3), (4.5) imply  $R  \in  C^{\infty}({\Bbb R}\setminus I_a,L(X))$ and,  with  $k\in \N_0= \N \cup \{0\}$,

(4.6) \qquad    $R^{(k)}(t) = k! (-i)^k R(t)^{k+1}$, $t \in {\Bbb R}\setminus I_a$,

(4.7) \qquad    $||R^{(k)}(t)|| \le  \eta_k |t|^{-(k+1)\delta}$, $t \in {\Bbb R}\setminus I_a$,
with  $\eta_k  : = k! \eta^{k+1} $, $\eta$ of (4.5).
\endproclaim

\demo{Proof} See [2, Appendix B, Corollary B.3, p. 463]. \P
\enddemo

\proclaim{Lemma 4.3} If  (4.1), (4.3)  hold and $F$ is a closed set $\st {\Bbb R}$ with $F\cap K =
\emptyset$, $K$ of (4.2), there exist  $H \in C^{\infty}({\Bbb R},L(X))$ and $M$ compact with  $K\st M\st {\Bbb R}\setminus F$,
so that  $H = R$  on  ${\Bbb R}\setminus M $ and $H(s)x \in D(A)$ for all $s \in{\Bbb R}$, $x \in X$.
\endproclaim

\demo{Proof} Lemma 3.5 with $U = {\Bbb R}\setminus F$ gives $\va\in\h ({\Bbb R})$ and $V$ open with $K\st V\st \overline {V}$ compact $\st U$,
$\va = 1$ on $V$, supp $\va \st {\Bbb R}\setminus F$.
Define   $H  : = R - \va R$  on  ${\Bbb R}\setminus K$,  $H  : = 0$  on $K $.
Then  $H = 0$ on the open $V$; if $t \not\in V$, then $t$ has positive distance to $K$,
so  $H = R - \va R$ on some  $(t-\e,t+\e)\st {\Bbb R}\setminus K$; so $H\in C^{\infty}
({\Bbb R},L(X))$ with Lemma 4.2. With  $M : =$ supp $\va \st {\Bbb R}\setminus F $, $K\st M$, one has
$H = R $ on  ${\Bbb R}\setminus M$.
$ R(\lambda,A) = (\lambda - A)^{-1}$ and the definition of  $(\lambda - A)^{-1}$ give
$  R(\lambda,A)x \in D(A)$ if $\lambda \in \cc\setminus \sigma(A)$, $x \in  X$; since $D(A)$ is linear,
  the definition of $H$ gives $H(s)X \st  D(A)$, $s \in{\Bbb R}$.  \P
\enddemo

\proclaim{Lemma 4.4} If (4.1), (4.3), (4.5) hold and $F$, $H$ are as in Lemma 4.3, then

(4.8) \qquad  $G(t)  : =  \frac{1}{2 \pi}  \int^{\infty}_{-{\infty}} H(s) e^{its}\, ds$

\noindent exists as an improper Riemann integral for all $t\in {\Bbb R}\setminus \{0\}$,
$H^{(k)}\in L^1({\Bbb R},L(X))$
 for $k\in\N$  and

(4.9) \qquad $G(t) = \frac{1}{2\pi} (i/t)^k \int_{{\Bbb R}} H^{(k)}(s) e^{its}\, ds$,  $k\in\N$, $t\in {\Bbb R}\setminus \{0\}$.

\noindent Moreover, $G\in L^1({\Bbb R}\setminus (-1,1),X)$ and is continuous at each $t\in {\Bbb R}\setminus \{0\}$.
\endproclaim
\demo{Proof} Partial integration yields for $t\not=0$, $S  < T$

(4.10)\qquad  $\int^T_S H(s) e^{its}\,ds = (1/it) (H(T)e^{itT} - H(S)e^{itS} -
             \int^T_S H'(s) e^{its}\, ds)$.

\noindent   With Lemma 4.3 there exists $b \in (0,\infty)$ with  $H = R$ on ${\Bbb R}\setminus (-b,b)$, (4.7)
gives $H(s) \to 0$ as  $|s| \to \infty$  and  $H^{(k)} \in  L^1({\Bbb R},L(X))$, $k\in\N$.
So the integral in (4.8) exists as $lim_{T\to \infty, S\to -\infty} \int^T_S H(s) e^{its}\, ds$,
(4.9) holds for  k=1. Induction gives (4.9) for  $k\in\N $.
Continuity of $G$ follows from (4.9), $H'\in L^1({\Bbb R},L(X))$ and the Lebesgue
dominated convergence theorem; $k=2$ and $H''\in L^1 ({\Bbb R},L(X))$ show  that
 $G \in L^1({\Bbb R}\setminus (-\e,\e),L(X))$  for any $\e >0$.
   \P \enddemo

\proclaim {Proposition 4.5} For $ H, G$ as in Lemma 4.4 one has  $G \in L^1({\Bbb R},L(X))$.
\endproclaim
\demo{Proof} By Lemma 4.4 it is enough to show integrability  of $G$ over $[-1,1]$.
By Lemma 4.3 there is $ b\in [max\{1,a\},\infty)$ with  $H(s)=R(s)$ if $|s|\ge b$,
$a$ of (4.4). Then

$2 \pi G(t) = \int^b_{-b} H(s) e^{its}\,ds + \lim_{T\to \infty}\int^T_b R(s)
      e^{its}\,ds   +  \lim_{S\to -\infty} \int^{-b}_S R(s) e^{its}\, ds $

  $ = :
       B_0(t) + B_+(t) + B_{-}(t) $ for  $t\in {\Bbb R}\setminus \{0\}$.

  $B_0 \in   L^1([-1,1],L(X))$ :  $B_0 (t)=\int^b_{-b} H(s) e^{its}\,ds$ is continuous on $[-1,1]$.

$B_+ \in   L^1([0,1],L(X))$: Set  $U(t) : =  \int^T_b R(s)
      e^{its}\,ds    = (1/t) \int^{tT}_{tb} R(u/t)e^{iu}\,du =
       (1/t)\int^b_{tb} R(u/t)e^{iu}\,du + (1/t)\int^{tT}_b R(u/t)e^{iu}\,du = :
          U_*(t)  +  U_T(t)$
 if $0< t\le 1$, $b < tT$.

 With (4.7) and $\frac{1}{2} < \delta < 1$ one has, independent of $T$,

 $ ||U_*(t)|| \le (1/t) \eta \int^b_{tb} (t/u)^{\delta} du \le
         ((\eta b^{1-\delta})/(1-\delta)) t^{\delta-1} = :  v_1 (t)$.

   Using  (4.6), $\frac{1}{2} < \delta < 1$ and partial integration,

   $||U_T(t)||  =  ||(1/(it))(R(T)e^{itT} - R(b/t)e^{ib} - \int^{tT}_b R'(u/t)
                     (1/t)  e^{iu}\, du )  ||  \le$

 $(1/t)(||R(T)|| + ||R(b/t)|| ) + (1/(t^2)) \int^{tT}_b ||R(u/t)||^2\, du $.

 $R(T)\to 0$  as $T\to\infty$ by (4.7),
  $(1/t) || R(b/t) || \le  \eta b^{-\delta} t^{\delta - 1}  = : v_2(t) $,

  $(1/(t^2)) \int^{tT}_b ||R(u/t)||^2 \,du \le(1/ (t^2)) \eta^2 t^{2\delta}
            \int^\infty_b u^{-2 \delta} du \le$

            $\eta^2 t^{2\delta -2}b^{1-2\delta}/(2\delta-1)  = : v_3(t)$.

 Together one gets for  $0<t\le 1$

 $ ||B_+(t)|| = ||U_*(t) + \lim_{T\to \infty} U_T(t)|| \le v_1(t)+v_2(t)+v_3(t)$
   with $(v_1+v_2+v_3) \in L^1([0,1],{\Bbb R})$ because $\delta > \frac{1}{2}$. Since $U = U_*
   + U_T$ is continuous on $[0,1]$, $B_+  = \lim_{T\to \infty} (U_* + U_T)$  is
Bochner measurable,  so  $B_+  \in  L^1([0,1],L(X)) $.

 $ B_{-} \in L^1[(0,1],L(X))$ follows from the above, since  $R(-t)$ also satisfies
   (4.5),  so  $G  \in   L^1([0,1],L(X))$.

   This for $H(-s)$ instead of $H(s)$  gives  $G  \in  L^1([-1,0],L(X))$,  and so
    $G  \in  L^1([-1,1],L(X)) $.
 \P
\enddemo

\proclaim{Corollary 4.6} For F, H, G as in Lemma 4.4 one has $\widehat {G} = H $  on ${\Bbb R}$.
\endproclaim
\demo{Proof} With $L^1 := L^1({\Bbb R},L(X))$, Lemma 4.3 and (4.7) one has ´$H^{(k)}
\in  L^1 \cap  C^{\infty}$ if $k\ge 1$, (4.9) gives

(4.11)\qquad $\widehat  {H'}(t)  =  2 i t\pi G(-t)$, $0\not =t \in {\Bbb R}$.

$G \in L^1$ of Proposition  4.5 and (4.9) for $k=3$  imply  $G_1 \in L^1$, $G_1(t) : =
t G(t)$, $t\not =0$; $G_1 \in C({\Bbb R},L(X))$ by (4.11), with $G_1(0)=\widehat  {H'}(0)/(-2 i\pi) =
0 $ with Lemma 4.3 and (4.7). With Proposition 3.7 for $f=H'$ one gets therefore

$ H'(s) = i \int_{{\Bbb R}} t G(-t) e^{ist}\,dt  = -i  \int_{{\Bbb R}} tG(t)e^{-ist}\,dt = -i\widehat{G_1}(s)$.
Since $G_1 \in L^1$, the dominated convergence theorem gives existence of
$\widehat{G}'$ on ${\Bbb R}$ and  $\widehat{G}'(s) = -i \widehat{G_1}(s)$, $s \in {\Bbb R}$, or  $\widehat{G}' = H'$ on  ${\Bbb R}$.
Since $H$ and $G$ vanish at infinity by (4.7) respectively  (4.9), $\widehat{G} = H$ on  ${\Bbb R}$.  \P
\enddemo

 \head {\S 5.  Existence of bounded uniformly continuous solutions of
         $u' = A u + \phi$    on  ${\Bbb R}$}\endhead

\proclaim{Lemma 5.1} Assume $D(A)$ linear $\st X$, $A : D(A) \to\,X$ linear, $x\in X$,
     $\lambda\in {\Bbb R}$ with  $i \lambda\in \cc\setminus \,\sigma(A)$, $e_{\lambda,x}(t) : =
    e^{i \lambda t} x$, $t\in {\Bbb R}$,  $F\in L^1({\Bbb R},L(X))$ with $\widehat{F}(\lambda) =
    R(i \lambda,A)$,

(5.1)\qquad  $v(t) : = (F*e_{\lambda,x})(t)  =  e^{i \lambda t}R(i \lambda,A)x$, $t\in {\Bbb R}$.

\noindent      Then v is a classical solution of

(5.2)\qquad     $v'  =  A v  +  e_{\lambda,x}$    on ${\Bbb R}$.
\endproclaim

\demo{Proof}  $R(i \lambda,A)x = (i \lambda - A)^{-1} x \in D(A)$ by definition of
$(i \lambda - A)^{-1}$, so $v(t)\in D(A)$, $t\in {\Bbb R} $. Since  $F\in L^1 ({\Bbb R},L(X))$  and $x \in X$ one has  $Fx \in L^1({\Bbb R},X)$  and
   $\widehat{(Fx)}  = (\widehat{F})x$ by [2, Proposition 1.1.6].  With $\widehat{F}(\la)= R(i\la,A)$
and Proposition 3.2 the  $v  =
F*e_{\lambda,x}$ is well defined, the following Bochner integrals all exist
and one has

   $v(t)=(F*e_{\lambda})(t) = e^{i\lambda t} \int_{{\Bbb R}} F(s)x  e^{- i \lambda s} ds =
   e^{i \lambda t}\widehat {(Fx)}(\lambda) = e^{i \lambda t}\widehat{F}(\lambda)x=$

 \qquad \qquad\qquad \qquad\,\,\,  $
  = e^{i \lambda t} R(i \lambda,A) x$,\qquad
  $t\in {\Bbb R}$, with $R(i \lambda,A)x \in D(A)$.

\noindent This implies  $v'(t) - A(v(t)) = e^{i \lambda t}(i \lambda
- A) R(i \lambda,A)x  =  e_{\lambda,x}(t)$, so  (5.2) holds.   \P
\enddemo

\proclaim{Theorem 5.2} Let  $T$ be a $C_0$-semigroup  on $X$ with generator $A$ satisfying (4.3)
and (4.4), $F$ closed  $\st {\Bbb R}\setminus K$  with $K = (1/i)(i{\Bbb R}) \cap \sigma(A))$. Then there exists $ G \in  L^1({\Bbb R},L(X))$ so that for any  $\phi\in L^{\infty}({\Bbb R},X)$ with  Beurling spectrum  $sp(\phi)\st F$ the  $u : = G*\phi$ is a mild solution of

(5.3) \qquad    $u'    =   A u  +   \phi$     on  ${\Bbb R}$

\noindent which is bounded and uniformly continuous on ${\Bbb R}$ with $sp (u) \st sp(\phi)\st F$.
\endproclaim

\demo{Proof} By Lemma 4.3  there exist $H\in \,C^{\infty}({\Bbb R},L(X))$ and $M$ compact with
$K\st M\st {\Bbb R}\setminus \,F$  and  $H = R$  of (4.2) on ${\Bbb R}\setminus \,M$. By Proposition 4.5 /Corollary 4.6 there exists  $G
\in L^1({\Bbb R},L(X))$  with  $\widehat{G} = H$ on ${\Bbb R}$, so  $\widehat{G} = R$  on ${\Bbb R}\setminus \,M $.
By Lemma 3.6 there exist  trigonometric polynomials $\Pi_n$ with $\Pi_n\to \phi$
almost everywhere on ${\Bbb R}$,  sup$_{n\in \N}||\Pi_n||_{\infty}  <  \infty$,  and all Fourier exponents
$\lambda$ of the $\Pi_n$  satisfy  $\lambda \in \,{\Bbb R}\setminus \,M$, so  $\widehat{G}(\lambda) = R(\lambda)$.
By Lemma 5.1 and linearity of $D(A)$ and $A$  the  $u_n : = G*\Pi_n $ are
classical and so mild solutions of (5.3) on ${\Bbb R}$ with $\Pi_n$ instead of $\phi$,
$u_n\in \,BUC({\Bbb R},X)$ by Proposition 3.2. By Lemma 3.1 the $u_n$ satisfy

(5.4)\qquad $u_n (t)= T(t-t_0)u_n(t_0)+\int^t_{t_0} T(t-s)\Pi_n(s)\, ds$, $t\ge  t_0$, $t_0 \in {\Bbb R}$.

\noindent If  $n \to \infty$, by Lebesgue's dominated convergence  theorem  $u_n \to
G*\phi$ on ${\Bbb R}$ with  $||u_n||_{\infty} \le  ||G||_{L^1} ||\Pi_n||_{\infty}$; since
sup $\{||T(s)|| : 0 \le s\le r\} < \infty $ for any $r\in \,{\Bbb R}_+ $ by the uniform boundedness theorem, one gets (5.4)  for $u = G*\phi$  and  $\phi$. Again Lemma 3.1 shows  that
$G*\phi$  is a mild solution of  (5.3) on ${\Bbb R}$, with $G*\phi\in \,BUC
({\Bbb R},X)$  by  Proposition 3.2 and $sp (G*\phi) \st sp(\phi)$ by Lemma 3.4.   \P
\enddemo

\proclaim{Remark 5.3 (a)}  The case  $\sigma  (A)  \cap  i{\Bbb R}  = \emptyset$, $F = {\Bbb R}$ is included,
      with  $H(t) = R(it,A)$, $t \in {\Bbb R}$,  in (4.8); then Theorem 5.2 can be applied
       with one $G$ for all $\phi \in L^{\infty}({\Bbb R},X)$.

     (b)  $\sigma  (A)  \cap  i{\Bbb R}  = \emptyset$ holds if $T$ admits exponential dichotomy
           [20, p. 409] : Then  $\sigma(T(1)) \cap  \{z \in \cc : |z| = 1 \} = \emptyset$ by
           [12, p. 191, 3.12 Proposition (b)], $ e^{\sigma(A)}  \st  \sigma(T(1))$ by [18,
            p. 45, (2.6)].
            A special case of this is $T$ exponentially stable [12, p. 186, 3.1 Definition
              (a), p. 188, 3.5 Proposition].
\endproclaim

\noindent \text{\bf{Example 5.4}}. $The\, assumption\,$ $\phi \in L^{\infty}({\Bbb R},X)$ $\,cannot\, be\, weakened \,to
         \,Stepanoff$

         \noindent $norm\,\,$  $||\phi||_{S^1}=\,\, $ sup$_{t\in {\Bbb R}} \int_t^{t+1} ||\phi (s)||\,ds   < \infty$,
       $\,\, not\, \,even\,\, for\, \,holomorphic\, and \,expo$-

       \noindent $nentially\, stable$ $T$ :

       Let $X = \cc$, $T(t) =$ multiplication by $e^{-t}$, $t \in {\Bbb R} $.
    Then the  generator of $T$ is  $A = -1$.
      Define  $\phi : = \sum_2^{\infty} \phi_n$  with $\phi_n = n$ on $ [n,n+\frac{1}{n}]$, else 0.
     Then  $||\phi||_{S^1} \le 2$,  $(i sp_C (\phi)) \cap \sigma (A) = \emptyset$ and
  $ u_{\phi}(t) : = \int_0^{\infty} T(s)\phi(t-s)\,ds$
 is a bounded mild solution  of  $ u' =Au+\phi$ on ${\Bbb R}$;  but no mild solution of this equation  on ${\Bbb R}$ or ${\Bbb R}_+$ is uniformly continuous.

   \demo{Proof} Indeed,   $ u_{\phi}(t)$    bounded follows as  in
     [4, p. 69], it is a mild solution with  $ u_{\phi}(t)= \int_{-\infty}^t
      T(t-s)\phi(s)\, ds$ and Lemma 3.1.  Then

       $ u_{\phi} (t_0+h) -   u_{\phi} (t_0) =(e^{-h} - 1) u_{\phi}(t_0) +  \int_{t_0} ^{t_0+h} T(t_0 +h-s)\phi (s)\,ds $,
     $ t_0 = n$, $h =\frac{1}{n}$

\noindent     and the boundedness of $u_{\phi}$ show that $u_{\phi}$ is not uniformly continuous on ${\Bbb R}_+$. \P
\enddemo

\proclaim {Corollary 5.5}  If $T$, $A$ are as in Theorem 5.2, $\phi\in \,L^{\infty}({\Bbb R},X)$ with
    $(i sp (\phi)) \cap \sigma (A) = \emptyset$,  then (5.3) has a mild solution  $u_{\phi}\in \,     BUC({\Bbb R},X)$.
\endproclaim

\proclaim {Corollary 5.6} If $T$ , $A$, $\phi$ are as in Corollary 5.5 and  sup$_{0<t< \infty}||T(t)||
         < \infty$, then for each  $x\in \,X$ the Cauchy problem  $u' = Au + \phi$ on
         ${\Bbb R}_+ $, $u(0) = x$ , has a unique solution, all these are in $BUC({\Bbb R}_+,X)$.
\endproclaim

\demo{Proof} [2, Proposition 3.1.16] and Corollary 5.5. \P
\enddemo
\proclaim {Corollary 5.7} If $T$, $A$, $\phi$ are as in Corollary 5.5, and $T$ is in addition a $C_0$-
   group, then for each $x \in X$ the equation (5.3) has a unique mild solution $u$
   on   ${\Bbb R}$ with $u(0) = x$; all these $u$ are $\in  BUC({\Bbb R},X)$ if  furthermore sup$_{t \in {\Bbb R}}
    ||T(t)||  <  \infty$.
\endproclaim
\demo{Proof} With Theorem 5.2 one can assume $\phi = 0$.

\noindent Uniqueness: (3.3) gives  $u(0)=T(0-(-n))u(-n) = T(n)u(-n)$,\,\,\, $u(-n) $

  $ = T(-n)u(0)$, $u(t) = T(t-(-n))T(-n)u(0) = T(t)u(0)$, $t > -n$, $n \in \N $.

\noindent   Existence:  $u(t) : = T(t)x $, $t \in{\Bbb R}$, gives  $u(t_0) = T(t_0)x$ or $x =T(-t_0)
   u(t_0)$,  so $u(t) = T(t)T(-t_0)u(t_0) = T(t-t_0)u(t_0)$, with Lemma 3.1
   $u$ is a mild solution   on ${\Bbb R}$.    \P
\enddemo

\proclaim {Remark 5.7} For $C_0$-semigroups backward uniqueness of the Cauchy
   problem for (5.3) holds if and only if all $T(t)$ are injective, $t \in{\Bbb R}_+$.
\endproclaim

\noindent So backward uniqueness is already false for $X = BUC({\Bbb R}_+,\cc)$, $T(t)f  =$
translate $f_t$,  $A = \frac{d}{ds}$.

\head{\S 6.  Existence  of generalized almost periodic solutions of equation (3.1) in the non-resonance case}\endhead

For $\A\st L^1_{loc}(\jj,X)$ where $\jj\in \{{\Bbb R},{\Bbb R}_+\}$, we define mean classes $\m\A$  by ([5, p. 120, Section 3 ])

(6.1)\qquad$\m\A:=\{f\in L^1_{loc}(\jj,X):   M_h f\in  \A ,\,\, h >0\}$, where

 \qquad \qquad \,\,\, $(M_h f) (t)=(1/h)\int_0^h f(t+s)\, ds$.

\noindent Usually, for example for $\A = AP$, $AA$, $VAA$, Stepanoff-, Besicovitch-,
         Eberlein weakly -, Levitan - almost periodic functions, recurrent
         functions one has  $\A  \st  \m\A\st \m^2\A\st \cdots$  with  the $\st$ in general strict (see [5, (3.8)], [7, (1.9)]).

  We denote by $\n$ any class of functions having the following properties:

(6.2)  $\n$ linear $\st L ^1_{loc}(\jj,\Bbb{X})\st X^{\jj}$.

(6.3) $(\phi_n)\st \n \cap  BUC$ and $\phi_n \to \psi $ uniformly on $\jj$ implies $\psi\in \n$.

(6.4) $\phi \in BUC ({\Bbb R},X)$, $\phi|\jj \in \n$ implies $\phi_a |\jj \in \n$ for $a\in {\Bbb R}$ ($\n$ BUC-invariant).

 (6.5) $ B\circ\phi|\jj \in \n$ for each $B\in L(\Bbb{X})$, $\phi\in BUC ({\Bbb R},X)$ with $ \phi|\jj\in  \n$.

\proclaim{Lemma 6.1}  If   $\n$ satisfies (6.2)-(6.5) and  $\phi\in L^{\infty}({\Bbb R},X)$ with $\phi|\jj\in \m\n$, $F\in L^1({\Bbb R},L(X))$ respectively $L^1 ({\Bbb R},\cc)$ then $F*\phi \in BUC({\Bbb R},X)$ with $(F*\phi)|\jj\in \n $.
\endproclaim

\demo{Proof} By Proposition 3.2 (i), $F*\phi$ exists and $\in BUC({\Bbb R},X)$.
To $F$ there is a sequence of $L(X)$-valued step-functions $H_n= \sum_{j=1}^{m_n} B_j\, \chi_{[\al_j,\beta_j)} $  with
$||F - H_n|| _{L^1}  \to 0$,  so  $||F*\phi  -  H_n*\phi||_{\infty} \le ||\phi||_{\infty} ||F - H_n|| _{L^1} \to 0 $ as $n\to \infty$. With (6.2), (6.3)
it is enough to show  $((B \chi_I)*\phi)|\jj \in \n$, $I = [\al,\beta)$,  for each $B\in L(X)$. With [2, Proposition 1.1.6]
one has $(B \chi_I)*\phi = B(\chi_I *\phi)$, with (6.5)  we have to show $\psi_{\al,\beta}|\jj\in \n$, where
$\psi_{\al,\beta} : = \int_{\al}^{\beta} \phi( \cdot - s) \,ds = \int_{-\beta}^{-\al}
\phi( \cdot + s)\,ds = \int _0 ^{- \al} \phi( \cdot + s) \,ds  -  \int_{-\beta} ^0 \phi( \cdot + s) \,ds $.
Now  $\phi|\jj \in\m\n $  gives  $\int_0^h \phi( \cdot + s) \,ds |\jj  \in \n$  if  $h>0$, (6.4) for $a = -h$
then  $\int_0^h \phi(\cdot - h + s)\, ds|\jj =  - \int_0^{-h} \phi( \cdot + s)\, ds |\jj  \in\n$, which
gives $\psi_{\al,\beta}|\jj \in\n$.

 The proof for  $F\in L^1 ({\Bbb R},\cc)$ is the same. \P
\enddemo

\proclaim {Examples 6.2} Examples of $\n$ satisfying (6.2)-(6.5) are,  for $\jj\in \{{\Bbb R},{\Bbb R}_+\}$, the  spaces of almost periodic  functions $AP=AP ({\Bbb R},X)$ [3,  2.1], Veech almost automorphic functions $VAA({\Bbb R},X)$, Bochner almost automorphic functions $AA({\Bbb R},X)$, Stepanoff $S^p$- almost  periodic functions, $1
\le p < \infty$ [5,p. 132],
          bounded Levitan almost periodic functions $LAP_b({\Bbb R},X)$ [15, Sec. 4, p. 53], [6, p. 430], linear subspaces with (6.4) of  bounded recurrent functions  of $REC_b({\Bbb R},X)=RC$ of [6, p. 427],  weakly almost periodic functions $\{f \in L^1_{loc}({\Bbb R},X) : y \circ f  \in
            AP({\Bbb R},\Bbb{C})$  for all $y \in$  dual $X^*\}$,  Eberlein almost  periodic functions $EAP(\jj,X)= EAP_0 (\jj,X)\oplus AP(\jj,X)$ [3, 2.3], asymptotic  almost periodic functions
           $AAP (\jj,X) = C_0 (\jj,X)\oplus AP(\jj,X)$, asymptotic almost automorphic functions $AAA(\jj,X)
           =  C_0 (\jj,X)\oplus AA(\jj,X)$, Zhang's (generalized) pseudo almost periodic functions  (G)PAP [27, p. 57, 67], pseudo almost automorphic functions  $PAA$ [24], [6,  Proposition 1.2, Examples 5.4/5.6].
\endproclaim

\proclaim{Theorem 6.3} Let $A$,  $\phi$ be as in Theorem 5.2, $\jj\in \{{\Bbb R},{\Bbb R}_+\}$ and $\phi|\jj \in \m\n$
with $\n$ satisfying  (6.2)-(6.5). Then there exists $ G \in  L^1({\Bbb R},L(X))$ so that  the  $u_{\phi} : = G*\phi$ is a mild  solution on ${\Bbb R}$  of (5.3)  and $ u_{\phi} \in  BUC({\Bbb R},X)$ with  $ u_{\phi}|\jj \in \n $.
\endproclaim
\demo{Proof} By Corollary 5.5, $u_{\phi}$ is a mild solution  of (5.3)    on ${\Bbb R}$ which is   bounded and uniformly continuous. Since $\phi|\jj \in \m\n $, by Lemma 6.1 $ u _{\phi} \in  BUC({\Bbb R},X) $  with  $ u_{\phi} |\jj \in \n$. \P
\enddemo

\proclaim{Remarks 6.4 (a)}  If T as in Theorem 5.2 admits exponential dichotomy, one has
    $i\, sp (\phi)  \cap \sigma(A) = \emptyset$ for any  $\phi \in L^{\infty} ({\Bbb R},X)$ (Remark 5.3(b)),
  then there exists $G \in  L^1({\Bbb R},L(X))$ so that the solution $G*\phi \in \n  \cap BUC({\Bbb R},X)$
   for all   $\phi \in L^{\infty} ({\Bbb R},X)$  with $\phi|\jj  \in  \m\n$; any of the $\n$  of
  Examples 6.2  can be used here.

(b)  The assumption $\phi|\jj \in \m\n$ in Theorem 6.3 cannot be replaced by  $\phi|\jj \in \n$ unless $\n$ satisfies $\n \st \m\n$. All the classes of Examples 6.2 satisfy this condition.
 However, for $\n$= the Banach space $\A_{g}= g\cdot AP$ with $g=e^{it^2}$,
$\phi = g$, $\jj={\Bbb R}$  and $T$, $A$ as in Example 5.4, all the assumptions of Theorem 6.3 are fulfilled
except $\phi \in \m\n$; though $g\in \n$,  no solution of  $u' = Au + g $ is in $\n$.

(c) In Theorem 6.3 for the case $\jj={\Bbb R}$ the assumption
`` $\phi  \in   \m \n $"
 can
   be generalized to
    ``$\phi \in   \m^n \n $ for some $n \in \N$",
  using
     Lemma 3.1 (ii). However, this is no real improvement, since for $\n$ with
     (6.2)-(6.5) and $L^{\infty} : = L^{\infty}({\Bbb R},X)$ one can show

 (6.6) \qquad $L^{\infty} \cap \m\n = L^{\infty} \cap \m^n\n = L^{\infty} \cap \h'_{L^{\infty}}$, $n \in \N$,

\noindent   where  $\h'_{\A} : = \{$ distribution $S \in \h'({\Bbb R},X): S*\va \in  \A$ for all $\va \in \h({\Bbb R})\}$
    of   [7, (1.7)].
\endproclaim

\proclaim{Example 6.5} $X = Y^n$, $Y$ complex Banach space, $A =$ complex $n \times n$ matrix,
$u$, $\phi$  $X$-valued in (3.1), $\phi\in L^{\infty}({\Bbb R},X)$. Then, if $i sp (\phi)$ contains no purely imaginary eigenvalue of $A$ and $\phi \in \m\n$, $\n$ with (6.2)-(6.5), then (5.3) has a   mild solution on ${\Bbb R}$ which belongs to $ \n \cap BUC$.
\endproclaim
Already this extends a well known result of Favard [13, p. 98-99] on the existence of almost periodic solutions for a $n^{th}$ order ordinary differential equation. See also [17, Example 3.4, p. 270- 271].

\demo{ Proof} The associated (semi)group $T $ is even entire,  $\theta = 1$ in (4.4). \P
\enddemo

 Another example would be a result on the almost periodicity of all
solutions of the inhomogeneous wave equation in the non-resonance
case [25, p. 179, 181 Th\'{e}or\`{e}me III.2.1], [2, Proposition 7.1.1], here one has a $C_0$-group, all
solutions of the homogeneous equation are almost periodic.

\indent School of Math. Sci., P.O. Box No. 28M, Monash University,
 Vic. 3800.

\indent E-mail "bolis.basit\@sci.monash.edu.au".

\indent Math. Seminar der  Univ. Kiel, Ludewig-Meyn-Str., 24098
Kiel, Deutschland.

\indent E-mail "guenzler\@math.uni-kiel.de".

\enddocument